\documentclass{amsart}
\usepackage{amsmath, amsthm}

\usepackage{amsbsy}
\usepackage{amssymb}
\def\be{\begin{equation}}
\def\ee{\end{equation}}

\def\bea{\begin{eqnarray}}
\def\eea{\end{eqnarray}}

\title{Euler's difference table and decomposition of tensor powers of adjoint representation of $A_n$ Lie algebra}

\author{A. M. Perelomov}

\address{Institute of Theoretical and Experimental Physics,\\  117259 Moscow, Russia}

\begin{document}
\maketitle

\begin{abstract}  
By using of Euler's difference table, we obtain simple explicit formula  for the decomposition of $k$-th tensor power of 
adjoint representation of $A_n$ Lie algebra  at $2 k \le{n+1}$.
\end{abstract}

\def\l{\lambda}
\def\R{{\Bbb R}}
\def\bchi{\boldsymbol{\chi}}
\def\ba{\boldsymbol{a}}
\def\ve{\varepsilon}
\def\a{\alpha}
\def\k{\kappa}
\def\l{\lambda}
\def\nn{\nonumber}
\def\ni{\noindent}
\def\bi{\bibitem}
\def\b{\hfil\break}

\section{Introduction}
Decomposition of tensor products of representations of simple complex Lie algebras is complicate problem.\footnote
{For the basic notations of Lie algebras see for example [Ov].}
For the Lie algebras of small rank $A_2, C_2, G_2$ and $A_3$ some results were obtained in  recent papers [FGP1] - [FGP4] 
(see also  references there). For the general case see, for example, papers [Ha],[BD] and references there.

In present note we consider the problem of decomposition of $k$-th tensor power of 
adjoint representation of $A_n$ Lie algebra  at stability domain $2 k \le{n+1}$ ,
where the results dont depend on $n$ [Ma].
Using Euler's difference table [Eu1]-[Eu2], (see also [Ri]) we obtain simple explicit formulas  
(19),\, (20) for such decomposition, which on best author knowlidge are new one.

\section{Euler's difference table}
In 1753 Euler with relation to the card game "Jeu de Recontre" investigated in details the permutations of $k$ numbers [1,2,...,k]
without fixed points [Eu1], [Eu2] (see also [Ri]). Such permutation has the name derangement and we denote the number 
of derangements as $d_k$. 

In order to find derangement numbers $d_k$ Euler constructed first the difference table from numbers $e_{k}^{j}$,\,
$j\le k$,\, that determined by basic recurrence relation
\be
e_{k}^{j}=e_{k}^{j+1}-e_{k-1}^{j},\, e_{j}^{j}  = j!\,,
\ee
and proved that
\be
e_{k}^{0} = d_k\,.
\ee

Note that $e_{k}^{j}$ is divided to $j!$ and it is useful to introduce higher derangement numbers
\be
d_{k}^{j}= \frac{e_{k}^{j}}{j!}\,.
\ee

For derangement numbers we have two basic recurrence formulae [9], [11]
\begin{eqnarray}
d_k &=& (k-1)(d_{k-1}+d_{k-2}),\quad d_0=1, d_1=0,\\
d_k &=& k d_{k-1}+(-1)^k,\quad d_0=1, 
\end{eqnarray}
and the generating function
\be \frac{e^{-x}}{1-x}=\sum _{k=0}^\infty \frac{d_k}{k!}\,x^k.
\ee

The first ten derangement numbers are
\be
\begin{array}{ccccccccccccc}
k&   \vert &0 &1 &2 &3 &4 &5 &6 &7 &8 &9 &10 \\
d_k& \vert &1 &0 &1 &2 &9 &44 &265 &1854 &14833 &133496 &1334961 \end{array} 
\ee

From (1) follows recurrence relation for numbers $d_n^k$
\be d_n^{k}=\frac1{k}\,(d_n^{k-1}+d_{n-1}^{k-1}).
\ee
Iterating it we obtain
\be d_n^k=\frac{1}{k!}\sum _{j=0}^k {k\choose j} d_{n-j}.    
\ee

Let us give also the generating function
\be \sum \frac{d_{n+k}^k}{n!}\,x^n=\frac{e^{-x}}{(1-x)^{k+1}} \ee
and the tables of numbers $e_n^k$ and $d_n^k$ 

\begin{equation}
\left[ 
\begin{array}{llllllllll}
1&0&0&0&0&0&0&0&0&0 \\
0&1&0&0&0&0&0&0&0&0 \\
1&1&2&0&0&0&0&0&0&0 \\
2&3&4&6&0&0&0&0&0&0 \\
9&11&14&18&24&0&0&0&0&0 \\
44&53&64&78&96&120&0&0&0&0 \\
265&309&362&426&504&600&720&0&0&0 \\
1854&2119&2428&2790&3216&3720&4320&5040&0&0 \\
14833&16687&18806&21234&24024&27240&30960&35280&40320&0 \\
133496&148329&165016&183822&205056&229080&256320&287280&322560&362880 \end{array}
\right]
\end{equation}

\begin{equation}
\left[ 
\begin{array}{llllllllll}
1&0&0&0&0&0&0&0&0&0 \\
0&1&0&0&0&0&0&0&0&0 \\
1&1&1&0&0&0&0&0&0&0 \\
2&3&2&1&0&0&0&0&0&0 \\
9&11&7&3&1&0&0&0&0&0 \\
44&53&32&13&4&1&0&0&0&0 \\
265&309&181&71&21&5&1&0&0&0 \\
1854&2119&1214&465&134&31&6&1&0&0 \\
14833&16687&9403&3539&1001&227&43&7&1&0 \\
133496&148329&82508&30637&8544&1909&356&57&8&1 \end{array}
\right]
\end{equation}

\section{Decomposition of tensor powers}

It is convenient to consider adjoint representation $ad$ of Lie algebra  $A_n $ as tensor 
\be X_i^j, \quad i,j=1,\ldots ,n+1,\ee
which satisfies the condition
\be 
\sum _{i=1}^{n+1} X_i^i=0.
\ee
The $k$th tensor power of $ad$ 
\be
X_k = X_{i_1}^{j_1}X_{i_2}^{j_2}\cdots X_{i_k}^{j_k}
\ee
has the decomposition
\be
 X_k\,=c_0^{k}Y_0+c_1^{k}Y_1+\cdots +c_k^{k}Y_k\,.
\ee
The quantities $Y_p\,, p=0,1,...,k$ decompose into irreducible representations of Lie algebra $A_n$ 
\be
Y_p = R_p + ...\,,\quad  R_p=X[p,0,...,0,p]\,.
\ee
The coefficients of terms denoted as ... are the Littlewood - Richardson coefficients , see [Ma].

The quantity $Y_p$  obtained by contractions in (15) $k-p$ upper indeces with $k-p$ lower indeces.
So we choose ordered subset of $p$ quantities $X_{i}^{j}$ from $k$ such quantities.
The number of such subsets is $\frac{1}{p!}{k\choose  p}$.
Then  we choose from $p$ quantities $X_i^j$  $l$ quantities of type $X_i^i$. This gives factor${p\choose  l}$.
The contraction on rest $k-l$ indeces gives the quantity $d_{k-l}$.
As result we have the formula
\be
c_p^k={k\choose p}\,\frac{1}{p!}\sum _{l=0}^{p}{p\choose l}\,d_{k-l}\,.
\ee
Taking into account (9) we obtain the main formula of present note

\be
c_j^{k} = {k\choose j} d_k^{j},
\ee
where ${k\choose j}$ are binomial coefficients,  $d_k^{j}$ are higher derangement numbers (3).

Note that coefficients $c_j^{k}$ satisfied the recurrence relation
\be c_{j+1}^{k}=\frac1{(j+1)^2}\left( (k-j)\,c_j^{k}+kc_j^{k-1} \right),\, c_{0}^{k} =d_k,
\ee
which is consequence of recurrence relation (8).

\section{Conclusion}
In conclusion we give the decomposition of first ten powers of adjoint representation of Lie algebra $A_n$ at $2 k \le{n+1}$.
\begin{eqnarray} 
X_1 &=& Y_1, \nn \\
X_2 &=& Y_0+2Y_1+Y_2,\nn \\
X_3 &=& 2Y_0+9Y_1+6Y_2+Y_3, \nn\\
X_4 &=& 9Y_0+44\,Y_1+42\,Y_2+12\,Y_3+Y_4, \nn\\
X_5 &=& 44\,Y_0+265\,Y_1+320\,Y_2+130\,Y_3+20\,Y_4+ Y_5, \\
X_6 &=& 265\,Y_0+1854\,Y_1+2715\,Y_2+1420\,Y_3+315\,Y_4\nn 
\end{eqnarray}
\[
+30\,Y_5+Y_6
\]
\begin{eqnarray}
X_7 &=& 1854\,Y_0+14833\,Y_1+25494\,Y_2+16275\,Y_3+4690\,Y_4\nn
\end{eqnarray}
\[
+651\,Y_5+42\,Y_6+Y_7
\]
\begin{eqnarray}
X_8 &=& 14833\,Y_0+ 133496\,Y_1+263284\,Y_2+198184\,Y_3+70070\,Y_4\nn
\end{eqnarray}
\[
+12712\,Y_5 +1204\,Y_6+56\,Y_7+\,Y_8
\]
\begin{eqnarray}
X_9 &=& 133496\,Y_0+ 1334961\,Y_1+2970288\,Y_2+2573508\,Y_3+1076544\,Y_4\nn
\end{eqnarray}
\[
+240534\,Y_5  +29904\,Y_6+2052\,Y_7+72\,Y_8+ Y_9
\]
\begin{eqnarray}
X_{10} &=& 1334961\,Y_0 +14684570\,Y_1+36377685\,Y_2+35636040\,Y_3 +17199210\,Y_4\nn
\end{eqnarray}
\[
  +4558428\,Y_5+699930\,Y_6
 +63240\,Y_7+3285\,Y_8+90\,Y_9+Y_{10} .
\]

%%%%%%%%%%%%%%%%%%%%%%%%%%%%%


\begin{thebibliography}{aaaa}

\bibitem[OV]{OV} A.L. Onishchik, E.B. Vinberg, 
{\em Lie Groups and Algebraic Groups}, Springer (1990)

\bibitem[FGP1]{OV} J. Fern\'andez N\'u\~{n}ez, W. Garc\'{\i}a Fuertes,  A.M. 
Perelomov, 
{\em J.Phys.A Math.Theor.}{\bf 47} 145202 (2014); arXiv: 1304.7203v2 [math-ph]

\bibitem[FGP2]{OV} J. Fern\'andez N\'u\~{n}ez, W. Garc\'{\i}a Fuertes,  A.M. 
Perelomov,
{\em J. Math. Phys.} {\bf 56} 041702 (2015); arXiv: 1405.2758v1[math-ph]

\bibitem[FGP3]{OV} J. Fern\'andez N\'u\~{n}ez, W. Garc\'{\i}a Fuertes,  A.M. 
Perelomov,
{\em J. Math. Phys.} {\bf 56} 091702 (2015); arXiv:1506.07815v1 [math-ph]

\bibitem[FGP4]{OV} J.. Fern\'andez N\'u\~{n}ez, W. Garc\'{\i}a Fuertes,  A.M. 
Perelomov,
{\em J. Nonlin. Math. Phys.} {\bf 25}, No.4, 618 (2018);arXiv:1705.03711v2 [math-ph]

\bibitem[Ha]{OV} P. Hanlon, {\em Adv. Math.} {\bf 56}, 238 (1985)

\bibitem[BD]{OV} G. Benkart, S. Doty,  arXiv:math/0108106 (2005)

\bibitem[Ma]{Ma} I. Macdonald {\em Symmetric Functions and Hall Polynomials}, 
2nd. Ed., Oxford Univ. Press (1995)

\bibitem[Eu1]{eu} L.Euler: {\em Calcul de la probabilite dans le jeu de recontre}, 
Mem.Acad. Sci. de Berlin {\bf 7}, 255-270  (1753)

\bibitem[Eu2]{Eu} L.Euler: {\em Solutio quaestionis curiosal ex doctrina combinationum}, 
Mem. Acad. Sci. de Berlin {\bf 3}, 57-64 (1779)

\bibitem[Ri]{Ri} J.Riordan : {\em An Introduction to Combinatorial Analysis}, N.Y., 
Wiley (1958)

\end{thebibliography}
\end{document}